\documentclass[11pt,a4paper]{article}
\usepackage{epsfig}
\usepackage{amssymb}
\usepackage{epic}
\usepackage{eepic}
\usepackage{amsmath}
\usepackage{enumerate}

\def\cT{{\mathcal T}}
\def\cI{{\mathcal I}}
\def\cJ{{\mathcal J}}
\def\cF{{\mathcal F}}

\def\P{{\mathbb P}}

\def\V{{\mathcal V}}

\newtheorem{theorem}{Theorem}[section]
\newtheorem{lemma}[theorem]{Lemma}
\newtheorem{corollary}[theorem]{Corollary}

\author{Erik I. Broman\footnote{Department of Mathematics, Uppsala University,
url: \texttt{http://www2.math.uu.se/\~{ }broman/}, e-mail:broman@math.uu.se
supported by the Swedish Research Council.}}

\date{\today}

\title{Stochastic Ordering of Infinite Geometric Galton-Watson Trees.}

\begin{document}
  \maketitle

  \begin{abstract}
We consider Galton-Watson trees with Geom$(p)$ offspring distribution. We
let $T_{\infty}(p)$ denote such a tree conditioned on being infinite. 
We prove that for any $1/2\leq p_1 <p_2 \leq 1$, there exists a coupling
between $T_{\infty}(p_1)$ and $T_{\infty}(p_2)$ such that 
$\P(T_{\infty}(p_1) \subseteq T_{\infty}(p_2))=1.$



  \end{abstract}

\section{Introduction} \label{secintro}
Consider a Galton-Watson tree $T(p)$ with some offspring distribution 
$\mu_p$ depending on the parameter $p.$ 
If $p_1<p_2$, then for many offspring distributions,
it is possible to couple the trees $T(p_1)$ and $T(p_2)$ 
such that $\P(T(p_1) \subseteq T(p_2))=1$. This is for example the case
with binomial, Poisson and negative binomial offspring distributions.
One can then ask whether this property is preserved under certain 
conditioning. For example, if we let $T_\infty(p)$
denote a sample of $T(p)$ conditioned on being infinite, is it the case 
that there exists a coupling of
$T_\infty(p_1)$ and  $T_\infty(p_2),$ where $ p_1<p_2,$ such that
\begin{equation}\label{eqn:conj}
\P(T_\infty(p_1) \subseteq T_\infty(p_2))=1? 
\end{equation}
Of course, if $p$ is such that the tree $T(p)$ is finite almost surely, 
then one needs to take care in defining $T_\infty(p).$ 

It is known (see Example 1.1 in \cite{B}) that there exists 
a parameterized offspring distribution
$\mu_p$ such that $T(p_1) \subseteq T(p_2)$ can be made to hold
almost surely, but for which \eqref{eqn:conj} fails.
However, \eqref{eqn:conj} is known to hold for the Poisson offspring
distribution (see Lyons, Peled and Schramm in \cite{LPS})
and for the Binomial$(d,p)$ offspring  distribution (see \cite{B})
when $d=2,3.$ It would of course be desirable to give a complete characterization
of the properties of $\mu_p$ needed to have \eqref{eqn:conj}. Unfortunately, this 
seems to be out of reach at the moment and the goal of this paper is 
to give another piece of the puzzle by studying 
the case of geometric offspring distribution
(i.e. when for some $p\in(0,1),$ the probability that an individual 
has $k$ offspring is $p^k(1-p)$).
Similar questions can be asked about percolation clusters on graphs. One of the 
most interesting cases is the corresponding problem for bond-percolation on 
${\mathbb Z}^d$;  see Open problem 4.3 of \cite{B}.

Another similar type of question was studied 
by Lusczak and Winkler in \cite{LW}. The authors considered trees $T_k,$
which are trees conditioned 
on being of a fixed size $k,$ and asked the natural question whether 
these conditioned trees could be coupled in increasing order, i.e. so that 
$\P(T_k \subset T_{k+1})=1$ (see Theorem \ref{thm:LW} below for a precise 
statement of this result). 
We will use their result as a key ingredient in the proof of our main result.

We proceed to give some definitions needed for the statements of the main 
results.
Let $o$ be the root, and consider the sets $\V_0=\{o\},$ 
$\V_n:=\{1,2,\ldots\}^n$ and
$\V=\{o\}\cup_{n=1}^\infty\{1,2,\ldots\}^n=\cup_{n=0}^\infty \V_n$. 
For any vertex $v=(v_1,\ldots,v_n)\in\V_n,$ with $n\geq 1,$ we define 
$v^-=(v_1,\ldots,v_{n-1})\in\V_{n-1},$ with the convention that if $n=1,$ 
then $v^-=\{o\}.$ We will think of $v^-$ as being the parent of $v.$
Furthermore if $v=(v_1,\ldots,v_n)\in\V_n,$ is such that $v_n \geq 2,$ we let 
$v'=(v_1,\ldots,v_n-1)\in \V_n$ and otherwise we let $v'$ be undefined. 
We think 
of $v'$ as being the 'youngest older sibling' of $v.$ 
Of course, such a vertex only exists 
if $v$ is not the first child of $v^-.$

For any two elements $u=(u_1,\ldots,u_k)$ and $ v=(v_1,\ldots,v_l)\in \V\setminus\{o\},$
we let 
$(u,v)=(u_1,\ldots,u_k,v_1,\ldots,v_l)$
denote the concatenation of $u$ and $v.$ For $i \geq 1,$ we will allow a slight
 abuse of notation, by writing
$(u,i)$ instead of $(u,(i))$ for $(u_1,\ldots,u_k,i).$
Furthermore, for any $u\in \V$
we let $(u,o)=(o,u)=u$ so that in particular, $(o,o)=o$.

Consider a set of vertices $V\subset \V.$ $V$ induces a set of
edges if we ''place'' an edge between $v$ and $v^-$ whenever $v,v^ -\in V.$
A tree $T$ is a then a connected graph consisting of a 
set of vertices $V(T)\subset \V$ and the induced set of edges $E(T).$ 
Since $V(T)$ determines the graph completely, we will 
leave $E(T)$ implicit. 
A {\em subtree} of a tree $T$ is defined to be a connected subgraph
of $T.$
For any $v\in \V$ and a tree $T,$ we let $T^v$ denote the subtree
(of $T$) with vertex set $V(T^v):=\{u\in V(T): u=(v,w) \textrm{ for some } 
w\in \V\}.$ Note that if $v\not \in V(T),$
we get that $V(T^v)=\emptyset.$ 
Informally, $T^v$ is simply the tree consisting of 
$v$ and the descendants of $v$ that belongs to $T.$
We also define $H(T^v):=\{w\in \V: (v,w)\in V(T^v)\},$ which is simply a 
shift of $T^v,$ mapping $v$ to $o.$ 
For $i\geq 1,$ we sometimes abuse notation and write $T^i$ instead of $T^{(i)}.$
We will let $|T|$ denote the number
of vertices of a tree $T$ and we call this the {\em size} of $T$.

We will let $\cF$ denote the set of trees $T$, with the following two properties.
Firstly, $o\in V(T), $ and secondly, if $v\in V(T)$ then $v'\in V(T)$ if $v'$ exists.
For $1\leq k<\infty,$ let $T_k$ be uniformly chosen among the trees $T\in \cF,$ 
such that $|T|=k.$ Let the distribution of $T_k$ be denoted by $\cT_k.$
We let $c_k$ denote the number of trees $T\in \cF$ such that $|T|=k.$
Thus, $c_1=1,c_2=1,c_3=2,c_4=5,\ldots$ The trees of $\cF$ are sometimes
called the rooted,
ordered trees, and it is well known that $(c_k)_{k \geq 1}$
are the Catalan numbers (see for instance \cite{Stanley}, Exercise 6.19).

The random trees with Geometric offspring distribution that we consider 
in this paper are formally defined as follows. For $p\in(0,1),$ 
let $(X_v)_{v \in \V}$
be an i.i.d. collection of random variables, indexed by $\V,$
such that $\P(X_v=k)=p^k(1-p)$
for $k=0,1,\ldots$. Then, let $V_0(T(p))=\{o\}$ and inductively 
let 
\begin{equation} \label{eqn:Vn}
V_{n+1}(T(p))=\bigcup_{v\in V_n(T(p))} \bigcup_{1\leq i\leq X_v}\{(v,i)\},
\end{equation}
for every $n\geq 0.$ 
We then define $T(p)$ to be the tree with vertex set
\[
V(T(p))=\bigcup_{n=0}^\infty V_n(T(p)).
\]
Observe that $V_n(T(p)) \subset \V_n$ so that  
this is the set of vertices of $T(p)$ at distance $n$ from the root.
Of course, if $X_v=0$ for some $v,$ the second 
union of \eqref{eqn:Vn} is over an empty set and so no descendants of $v$ belongs to 
$T(p).$
For a similar reason, if for some $n,$ $V_n(T(p))=\emptyset,$ it follows that
$V_{n+1}(T(p))=\emptyset$. 
We also note that $T(p)\in \cF.$

It is easy to see from the construction, that for any tree $T\in \cF$ 
of size $k,$
$\P(T(p)=T)=p^{k-1}(1-p)^k,$ and therefore
\begin{equation} \label{eqn:probTk}
\P(|T(p)|=k)=c_{k}p^{k-1}(1-p)^k.
\end{equation}
It follows that the distribution of $T(p),$ conditioned on the event 
that $|T(p)|=k,$ 
is $\cT_k$ (in particular, it is independent of $p$). Sometimes, it will be 
convenient to think of the empty set as the tree of size $0,$
and then we will use the notation $T_0=\emptyset.$
For every $1\leq k \leq \infty,$ we define $\eta_k(p):=\P(|T(p)|=k).$ 
Note that we allow for $k=\infty,$ and that $\sum_{k=1}^\infty\eta_k(p)+\eta_\infty(p)=1,$
since we always start with the root so that $\P(|T(p)|=0)=0.$

For two trees $T,T'$
we say that $T \subseteq T'$ if $V(T) \subseteq V(T').$ 
Note that it follows 
by definition that $E(T) \subseteq E(T').$  
The following theorem is due to \cite{LW}. 
It is not explicitly stated, but as they point out
(p. 427), it follows from their 
argument.
\begin{theorem}[Luczak, Winkler] \label{thm:LW}
There exists a coupling of $(T_k)_{k\geq 1}$
(where $T_k \sim \cT_k$ for every $1\leq k<\infty$)
such that
\begin{equation} \label{eqn:1}
\P(T_1 \subset T_2 \subset \cdots )=1.
\end{equation}
\end{theorem}
{\bf Remark:} 
It is proved in \cite{Janson} that 
equation (\ref{eqn:1}) does not hold for 
general offspring distributions.

\medskip

It is well known, that for $p > 1/2$, $\P(|T(p)|=\infty)>0$. For such $p>1/2,$
let $T_{\infty}(p)$ denote a random tree whose distribution equals 
that of $T(p),$ 
conditioned on the event $|T(p)|=\infty.$ Furthermore, with $(T_k)_{k\geq 1}$
as in Theorem \ref{thm:LW}, we define 
\begin{equation} \label{eqn:Tinftydef1/2}
T_\infty(1/2):=\bigcup_{k=1}^\infty T_k.
\end{equation}
For any $p\geq 1/2,$ we let $\cT_\infty(p)$ denote the distribution of
$T_\infty(p).$
Since $(T_k)_{k\geq 1}$ does not depend on the parameter $p$, it may 
seem strange to use the notation $T_\infty(1/2)$ (i.e. $p=1/2$). 
However, in Section \ref{sec:prel}, we will describe a way to sample
trees with distribution $\cT_\infty(p),$ and then we will see that 
$\bigcup_{k=1}^\infty T_k$ naturally corresponds to the critical case $p=1/2.$

We can now state our main theorem.
\begin{theorem} \label{thm:main}
For any $1/2\leq p_1<p_2\leq1,$ there exists a coupling of 
$T_{\infty}(p_1)$ and $T_{\infty}(p_2)$
(where $T_{\infty}(p_1)\sim \cT_{\infty}(p_1)$
and $T_{\infty}(p_2)\sim \cT_{\infty}(p_2)$) such that
\[
\P(T_{\infty}(p_1) \subseteq T_{\infty}(p_2))=1.
\]
\end{theorem}
{\bf Remark:} We prove the theorem by giving an explicit construction of the 
coupling.

As mentioned before, the corresponding result for  
Galton-Watson trees 
with Poisson offspring distributions was proved in \cite{LPS}, 
while the corresponding 
result for ${\rm Bin}(d,p)$ offspring distribution was proved 
in \cite{B} for $d=2,3.$

\medskip

We have the following corollary.
\begin{corollary} \label{corr:}
For any $k$ and $p \geq 1/2,$ there exists a coupling of 
$T_k$ and $T_{\infty}(p),$ 
(where $T_k\sim \cT_k$ and $T_{\infty}(p)\sim \cT_{\infty}(p)$ )
such that
\[
\P(T_k \subset T_{\infty}(p))=1.
\]
\end{corollary}
{\bf Remark:} This is an immediate consequence of Theorem \ref{thm:LW}, 
\eqref{eqn:Tinftydef1/2} and Theorem \ref{thm:main}.

\medskip

We will end this section by briefly discussing the difficulties involved in 
proving Theorem \ref{thm:main}. In order to construct a coupling satisfying 
the statement of Theorem \ref{thm:main}, it is useful to be able to 
generate (or sample) trees  with distribution $\cT_\infty(p).$ A natural way to
do this is to use a sequential procedure using conditional probabilities as 
follows (where we will use the easily established result, 
see Lemma \ref{lemmaTpinfinite}, that $\eta_\infty(p)=(2p-1)/p$
for $p>1/2$).
Informally, we construct $\bar T (p)\sim \cT_\infty(p),$ 
by in the first step letting $\bar T^1 (p)$
be infinite with probability 
$\P(|T^1(p)|=\infty | |T(p)|=\infty)=p\eta_\infty(p)/\eta_\infty(p)=p,$
and of size $k<\infty$ with probability 
$\P(|T^1(p)|=k| |T(p)|=\infty)=p\eta_k(p).$ In the second step,
if $|\bar T^1(p)|=\infty,$ then we let $\bar T^2(p)$ be infinite with 
probability $\P(|T^2(p)|=\infty | 
|T^1(p)|=|T(p)|=\infty)=p\eta_\infty(p),$
while if $|\bar T^1(p)|=k<\infty,$ we let $\bar T^2(p)$ be infinite with 
probability $\P(|T^2(p)|=\infty | 
|T^1(p)|=k, |T(p)|=\infty)=\cdots=p$. 
Continuing in this way until one finds a subtree which is infinite, and 
later a subtree which is of size $0$ (which marks the end of the procedure),
one can produce a tree $\bar T(p)\sim \cT_\infty(p).$

This is one of the most natural ways of sampling a tree with distribution $\cT_\infty(p).$
However, it is not possible to use this procedure to construct a coupling proving 
Theorem \ref{thm:main}.
This can be seen by first observing that any coupling 
of $T_\infty(p_1)$ and $T_\infty(p_2)$ satisfying 
$\P(T_\infty(p_1) \subseteq T_\infty(p_2))=1$ must certainly satisfy 
$|T^n_\infty(p_1)| \leq |T^n_\infty(p_2)|$ for every $n.$ 
With positive probability,
we could have that $|\bar T^1(p_1)|<\infty$ while $|\bar T^1(p_2)|=\infty.$
Then, the conditional 
probability that $|\bar T^2(p_1)|=\infty$ is $p_1,$ while the conditional 
probability that $|\bar T^2(p_2)|=\infty$ is $p_2\eta_\infty(p_2).$
Of course, 
for some choices of $p_1,p_2$ we can have that $p_1>p_2\eta_\infty(p_2),$
and so we do not get that $|T^2(p_1)| \leq |T^2(p_2)|$ with probability 
1.

In a second attempt, one might try to remedy the problem of our first attempt
by first determining which subtrees 
$\bar T^1(p),\bar T^2(p),\ldots$ should be infinite, and which should be finite.
Then, one could proceed by coupling 
these subtrees so that if $\bar T^i(p_1)$ is infinite than so is 
$\bar T^i(p_2).$ 
However, one will then find that with positive probability, both
$\bar T^1(p_1)$ and $\bar T^1(p_2)$ are finite, and then one 
would have to let $|\bar T^1(p_i)|=k$ with probability
\begin{equation} \label{eqn:Tcondfinite}
\P(|T(p_i)|=k | |T(p_i)|<\infty)=\frac{\eta_k(p_i)}{1-\eta_\infty(p_i)}
=c_kp_i^{k-1}(1-p_i)^k \frac{p_i}{1-p_i}=\eta_k(1-p_i).
\end{equation}
Because of 
\eqref{eqn:Tcondfinite}, we see that conditioned on $\bar T^1(p_1)$ and 
$\bar T^1(p_2)$ both being finite, 
we cannot have that $|T^1(p_1)| \leq |T^1(p_2)|$ with probability 1. In fact, a canonical 
coupling will result in $|T^1(p_2)| \leq |T^1(p_1)|.$

Therefore, the key to proving Theorem \ref{thm:main} is to find a way to generate
a tree with distribution $\cT_\infty(p)$ which has the desired monotonicity 
properties and does not fall into any of the traps described above. 
In Section \ref{sec:prel} we give this procedure along with 
some preliminary results. We then use this in Section \ref{sec:mainproof}
to prove Theorem \ref{thm:main}.

\section{Generating a tree with distribution $\cT_\infty(p)$} \label{sec:prel}

We start this section by proving the following easy lemma, 
already used in the introduction.
\begin{lemma} \label{lemmaTpinfinite}
If $p\leq1/2,$ then $\eta_{\infty}(p)=0$ and if $p>1/2,$ then
$\eta_{\infty}(p)=\frac{2p-1}{p}.$
\end{lemma}
{\bf Proof.} We have that
\begin{eqnarray*}
\lefteqn{1-\eta_\infty(p)=\sum_{k=0}^\infty \P(X_o=k)(1-\eta_\infty(p))^k}\\
& & =\sum_{k=0}^\infty p^k(1-p)(1-\eta_\infty(p))^k
=\frac{1-p}{1-p(1-\eta_\infty(p))}.
\end{eqnarray*}
Solving for $\eta_\infty(p),$ 
we have two solutions, $\eta_{\infty}(p)=0$ and
$\eta_{\infty}(p)=\frac{2p-1}{p}.$ It is well known that the tree is
supercritical iff $p>1/2$ from which the result follows.
\fbox{}\\

\medskip

We proceed by describing a procedure that will generate a tree $\tilde{T}(p).$
In Lemma \ref{lemma:welldefined} we show that this procedure is well defined.
Lemmas \ref{lemma:tildeTdistr1} and \ref{lemma:tildeTdistr2} will then prove
that $\tilde{T}(p)\sim \cT_\infty(p)$ for every $p\geq 1/2,$ while Lemma 
\ref{lemma:mon} will provide us with the crucial monotonicity 
properties used to prove our main results. 

However, before we give any details, we will explain some heuristics
of the construction and the coupling.
Consider therefore a tree with distribution $\cT_\infty(1/2).$ It is 
well known (see for instance \cite{Janson2} sections 5 and 7)
that such a tree will consist of one single infinite  path
to which there are smaller trees attached. In fact, it is
possible to prove (and indeed we do this in Lemma \ref{lemma:tildeTdistr1}) 
that such a tree can be generated in the following informally described 
way. 
Start with the root $o,$ and sample $X_1,X_2$ independently from a 
Geometric distribution with parameter $1/2.$ Then, attach independent 
trees with distribution 
$\cT(1/2)$ to the vertices $(1),\ldots,(X_1)$. Proceed by 
attaching an independent tree with distribution $\cT_\infty(1/2)$ 
to $(X_1+1)$
and end the construction by attaching independent trees with distribution 
$\cT(1/2)$ to the vertices $(X_1+2),\ldots,(X_1+X_2+1)$.
This implies that the probability that the infinite subtree can be found
at position $l$ (i.e. belongs to the lth child of the root) 
equals $\P(X_1=l-1)=2^{-l}.$

In order to construct the coupling of $\tilde{T}(1/2)\sim \cT(1/2)$ 
and $\tilde{T}(p)\sim \cT(p)$ for 
$1/2<p$ in the desired way, we start by finding the infinite subtree 
of $\tilde{T}(1/2)$ as above. For this, we simply use a random variable 
$X$ such that $\P(X=l)=2^{-l}.$ Given this number $l,$ 
the above paragraph 
tells us that for $\tilde{T}(1/2)$, the subtrees 
$\tilde{T}^1(1/2),\ldots , \tilde{T}^{l-1}(1/2)$ should be i.i.d. $\cT(1/2).$
We then require the corresponding subtree $\tilde{T}^l(p)$
to be infinite also for $\tilde{T}(p).$ 
The conditional distributions of $\tilde{T}^1(p),\ldots , \tilde{T}^{l-1}(p)$ are slightly
more complicated, but the key is that they will be independent, which facilitates a 
coupling of $\tilde{T}^i(1/2)$ and  $\tilde{T}^i(p)$ for $1\leq i \leq l-1$
such that $\tilde{T}^i(1/2) \subset \tilde{T}^i(p).$ Of course, it is not apriori
clear that it is possible to achieve the coupling $\tilde{T}^i(1/2) \subset \tilde{T}^i(p)$
but this is addressed in Lemma \ref{lemma:mon}.
We note that if our construction would instead condition on the smallest 
number $k$ such that 
$\tilde{T}^k(p)$ is infinite, then the conditional distributions of 
$\tilde{T}^1(p),\ldots \tilde{T}^{k-1}(p)$ would be $\cT(1-p)$ as described
in (\ref{eqn:Tcondfinite}). Thus, just as described 
in the Introduction, our coupling would fail. 

Of course, we also have to address the subtrees $\tilde{T}^{l+1}(1/2),\ldots$
and $\tilde{T}^{l+1}(p),\ldots$ which is done in a similar, albeit slightly more complicated
(as the conditional distributions of $\tilde{T}^{l+1}(p),\ldots$ are dependent of 
each other and $\tilde{T}^{1}(p),\ldots,\tilde{T}^{l}(p)$) way, to ensure that the end 
result $\tilde{T}(1/2)$ ($\tilde{T}(p)$) 
indeed has distribution $\cT_\infty(1/2)$ ($\cT_\infty(p)$). 
Finally, it will turn out that the construction
will not only work to couple $\tilde{T}(1/2)$ and $\tilde{T}(p),$
but will in fact work for any $1/2\leq p_1<p_2\leq 1.$

Below, $X$ will (as discussed) determine the position
of the single infinite subtree for $\tilde{T}(1/2)$, $L_m(p)$ will determine 
the size of the subtree $\tilde{T}^m(p)$ while $m_0-1$ will be the 
total number of subtrees.

\medskip

We now turn to the formal description.
We let $X$ be such that $\P(X=l)=2^{-l}$ for every $l\geq 1$,
and $(U_i)_{i \geq 1}$ be an i.i.d. sequence of $U[0,1]$ random
variables  which is also independent of $X.$ For $p\geq 1/2,$ we then define the 
random variables $L_1(p),L_2(p),\ldots$ 
through the following procedure

\begin{enumerate}
\item For every $m<X$ and $1\leq k < \infty,$ we let $L_m(p)=k$ if
\begin{equation} \label{eqn:tildeT1}
\sum_{l=1}^{k-1}2p\eta_l(p)\leq U_m <\sum_{l=1}^{k}2p\eta_l(p),
\end{equation}
while if $U_m >\sum_{l=1}^{\infty}2p\eta_l(p)$ we let $L_m(p)=\infty.$

\item We let $L_X(p)=\infty$.

\item In order to define $L_{m}(p)$ for every $m>X,$ we proceed sequentially. 
Assume therefore that $L_{m-1}(p)$ has been determined. 
We then define $n_\infty(p,m):=|\{i\leq m-1:L_i(p)=\infty\}|$ and let 
$L_m(p)=0$ if 
\begin{equation}\label{eqn:tildeT2}
U_m < \frac{2^{n_\infty(p,m)}p(1-p)}{(2^{n_\infty(p,m)}-2)p+1},
\end{equation}
and for $1\leq k < \infty,$ we let $L_m(p)=k$ if 
\begin{equation}\label{eqn:tildeT3}
\frac{2^{n_\infty(p,m)}p(1-p)}{(2^{n_\infty(p,m)}-2)p+1}+
\sum_{l=1}^{k-1}p\eta_l(p)\leq U_m <
\frac{2^{n_\infty(p,m)}p(1-p)}{(2^{n_\infty(p,m)}-2)p+1}+
\sum_{l=1}^{k}p\eta_l(p),
\end{equation}
and otherwise we let $L_m(p)=\infty.$
\end{enumerate}
For our construction of $\tilde{T}(p),$ we will only use the $L_m(p)$
such that $m<m_0:=\min\{m:L_m(p)=0\}.$ 
We have the following lemma.
\begin{lemma} \label{lemma:welldefined}
By the above construction, for any $m<l,$ we have that 
$\P(L_m(p)=\infty |X=l)=p\eta_\infty(p).$ Furthermore, we have 
that for $m>l,$ 
\[
\P(L_m(p)=\infty| X=l,n_\infty(p,m)=n )=p\eta_\infty(p)
\frac{(2^{n}-1)p}{(2^{n}-2)p+1}.
\]
\end{lemma}
{\bf Proof.} Using that $\eta_\infty(p)=(2p-1)/p,$ we have that 
\[
p\eta_\infty(p)+\sum_{k=1}^\infty 2p\eta_k(p)
=2p-1+2p(1-\eta_\infty(p))=2p-1+2(1-p)=1,
\]
and furthermore, for any $n\geq 1,$
\begin{eqnarray*}
\lefteqn{\frac{2^{n}p(1-p)}{(2^n-2)p+1}
+\frac{(2^n-1)p}{(2^n-2)p+1}p\eta_\infty(p)+\sum_{k=1}^\infty p\eta_k(p)}\\
& & =\frac{2^{n}p(1-p)}{(2^n-2)p+1}
+\frac{(2^n-1)p}{(2^n-2)p+1}(2p-1)+p(1-\eta_\infty(p))\\
& & =\frac{p^2(2^{n+1}-2-2^n)+p(-2^n+1+2^n)}{(2^n-2)p+1}+1-p\\
& & =p\frac{p(2^{n}-2)+1}{(2^n-2)p+1}+1-p=1,
\end{eqnarray*}
proving the lemma.
\fbox{}\\
{\bf Remark:} Note that this also proves that the procedure determining 
$(L_m(p))_{m\geq 1}$ is well defined, 
since all relevant (conditional) probabilities are 
positive and sum to one.

\medskip

We now turn to the construction of $\tilde{T}(p)$. 
\begin{enumerate}
\item Let $X,(L_m(p))_{m \geq 1}$ be defined as above.

\item Let $(T_{\infty,m}(p))_{m\geq 1}$ be an i.i.d. sequence which is independent of
the other random variables in this list, and where 
$T_{\infty,m}(p)\sim \cT_\infty(p)$ for every $m\geq 1.$

\item Let for every $m\geq 1,$ $(T_{k,m})_{k\geq 1},$ be a sequence of 
random variables such that $T_{k,m}\sim \cT_k$ for every $k\geq 1.$
Furthermore, let $(T_{k,m})_{k\geq 1}$ be independent for different $m$
and independent of the other random variables in this list.

\end{enumerate}

Recall that $m_0=\min\{m:L_m(p)=0\},$ and define the tree $\tilde{T}(p)$
by
\[
V(\tilde{T}(p))=\{o\}\bigcup_{m=1}^{m_0-1} \bigcup_{v\in T_{L_m(p),m}}\{(m,v)\}.
\]
Here, we abuse the notation somewhat in that we write $T_{L_m(p),m}$
instead of $T_{L_m(p),m}(p)$ when $L_m(p)=\infty.$
Informally, the tree $\tilde{T}(p)$ is constructed by starting with a root
and then attaching trees of size $L_m(p)$ at the vertex $\{(m)\}$ for every 
$m<m_0.$
One key property of the tree $\tilde{T}(p)$ is that 
$\tilde{T}(p) \sim \cT_\infty(p),$ which we prove next. Since the proofs
in the cases $p=1/2$ and $p>1/2$ are completely different, we split
the result into two lemmas. We point out that Lemma \ref{lemma:tildeTdistr1}
can be shortened by using the technique of size-biased Galton-Watson trees
at criticality (see \cite{Janson2}, Chapter 7, and the references within). 
However, in order to keep the paper self-contained, we give a proof from 
first principles.

\begin{lemma} \label{lemma:tildeTdistr1}
We have that $\tilde{T}(1/2) \sim \cT_\infty(1/2).$
\end{lemma}
{\bf Proof.}
Recall from \eqref{eqn:Tinftydef1/2} the definition of $T_\infty(1/2).$
We have from this 
that for any $1\leq l\leq m,$ $1\leq k_1,\ldots,k_{l-1},k_{l+1},k_m<\infty,$
and with $k_{l}=k_l(n)=n-(k_1+\cdots +k_{l-1}+k_{l+1}+\cdots+k_m)-1,$ 
\begin{eqnarray*}
\lefteqn{\P(|T_\infty^1(1/2)|=k_1,\ldots,|T_\infty^{l-1}(1/2)|=k_{l-1},
|T_\infty^{l}(1/2)|=\infty, }\\
& & \hspace{20mm} |T_\infty^{l+1}(1/2)|=k_{l+1},\ldots,
|T_\infty^{m}(1/2)|=k_{m},|T_\infty^{m+1}(1/2)|=0)\\
& & =\lim_{n \to \infty}\P(|T_n^1|=k_1,\ldots,
|T_n^{m}|=k_{m},|T_n^{m+1}|=0) \\
& & =\lim_{n \to \infty}\P(|T^1(p)|=k_1,\ldots,
|T^{m}(p)|=k_{m},|T^{m+1}(p)|=0 | |T(p)|=n) \\
& & =\lim_{n \to \infty}\frac{\prod_{i=1}^m p \eta_{k_i}(p)}{\eta_n(p)}(1-p)
=\lim_{n \to \infty}
\frac{p^{m+k_1-1+\cdots+k_m-1}(1-p)^{k_1+\cdots+k_m}}{p^{n-1}(1-p)^n}
\frac{\prod_{i=1}^mc_{k_i}}{c_{n}}(1-p)\\
& & =\lim_{n \to \infty}\frac{\prod_{i=1}^mc_{k_i}}{c_{n}}
=\left(\prod_{1\leq i\leq m: i\neq l} c_{k_i}\right)
\lim_{n \to \infty}\frac{c_{n-k-1}}{c_{n}},
\end{eqnarray*}
where we use \eqref{eqn:probTk} and $k:=k_1+\cdots +k_{l-1}+k_{l+1}+\cdots+k_m.$
Observe that the choice of $p\in(0,1)$ is irrelevant (as discussed in the introduction),
and that if $|T(p)|=n,$ then $\sum_{j=1}^\infty|T^j(p)|=n-1$
which explains the definition of $k_l.$
Furthermore, as in Lemma 2.1 of \cite{B}, we have that 
$\lim_{n \to \infty}c_{n-1}/c_n=1/4,$ so we conclude that 
\begin{eqnarray*}
\lefteqn{\P(|T_\infty^1(1/2)|=k_1,\ldots,|T_\infty^{l-1}(1/2)|=k_{l-1},
|T_\infty^{l}(1/2)|=\infty, }\\
& & \hspace{20mm} |T_\infty^{l+1}(1/2)|=k_{l+1},\ldots,
|T_\infty^{m}(1/2)|=k_{m},|T_\infty^{m+1}(1/2)|=0)\\
& &=\frac{1}{4^{k+1}}\left(\prod_{1\leq i\leq m: i\neq l} c_{k_i}\right)
=\frac{1}{4}\prod_{1\leq i\leq m: i\neq l} \frac{c_{k_i}}{2^{2k_i}}
=\frac{1}{4}\prod_{1\leq i\leq m: i\neq l}\frac{1}{2}\eta_{k_i}(1/2).
\end{eqnarray*}

We continue by noting that 
\begin{eqnarray*}
\lefteqn{\P(|\tilde{T}^1(1/2)|=k_1,\ldots,|\tilde{T}^{l-1}(1/2)|=k_{l-1},
|\tilde{T}^{l}(1/2)|=\infty, }\\
& & \hspace{20mm} |\tilde{T}^{l+1}(1/2)|=k_{l+1},\ldots, 
|\tilde{T}^{m}(1/2)|=k_{m}, |\tilde{T}^{m+1}(1/2)|=0) \\
& & =\P(L_1(1/2)=k_1,\ldots,L_{l-1}(1/2)=k_{l-1},L_l(1/2)=\infty,\\
& & \hspace{20mm} 
L_{l+1}(1/2)=k_{l+1},\ldots,L_{m}(1/2)=k_{m},L_{m+1}(1/2)=0 |X=l)\P(X=l)\\
& & =\left(\prod_{i=1}^{l-1}\eta_{k_i}(1/2)\right) 
\left(\prod_{i=l+1}^{m}\frac{1}{2}\eta_{k_i}(1/2)\right)\frac{1}{2}2^{-l}
=\frac{1}{4}\prod_{1\leq i\leq m: i\neq l}\frac{1}{2}\eta_{k_i}(1/2).
\end{eqnarray*}

This establishes that the joint distribution of the {\em sizes} of 
$T^1_\infty(1/2),T^2_\infty(1/2),\ldots$ is the same as that of
$\tilde{T}^1(1/2),\tilde{T}^2(1/2),\ldots$. Furthermore, it 
is easy to see that conditioned on $|T_\infty^i(1/2)|=k_i$, the 
distribution of $H(T_\infty^i(1/2))$ is $\cT_{k_i}.$ By construction 
of $\tilde{T}(1/2)$, the same holds for $\tilde{T}^i(1/2)$. 
Furthermore, by \eqref{eqn:Tinftydef1/2}, $\cT_{\infty}(1/2)$ is the limiting
distribution of $T_n\sim \cT_n$ when letting $n \to \infty.$ Therefore, 
conditioned on $|T_\infty^l(1/2)|=\infty$, we have that 
$H(T_\infty^l(1/2)) \sim \cT_{\infty}(1/2).$ By construction 
of $\tilde{T}(1/2)$, it holds that also 
$H(\tilde{T}^l(1/2))\sim \cT_{\infty}(1/2).$
We conclude that $T_\infty(1/2)$
and $\tilde{T}(1/2)$ must have the same distribution.
\fbox{}\\

\begin{lemma} \label{lemma:tildeTdistr2}
For every $p>1/2,$ we have that $\tilde{T}(p) \sim \cT_\infty(p).$
\end{lemma}
{\bf Proof.}
Let $\cI,\cJ\subset \{1,2,\ldots\}$
be such that $\cI\neq \emptyset,$ $\cI\cap \cJ =\emptyset,$ 
and $\cI\cup \cJ=\{1,2,\ldots,I+J\}$ where 
$I=|\cI|$ and $J=|\cJ|.$ Informally, $\cI=\{i_1,\ldots,i_I\}$ 
will be the set of children with an 
infinite number of descendants, while $\cJ=\{j_1,\ldots,j_J\}$ 
will be the set of children with a 
finite number of descendants. Here, we have ordered the elements 
of $\cI,\cJ$ so that $i_1<i_2<\cdots<i_I$ and
$j_1<j_2<\cdots<j_J.$ Using that $T_\infty(p)$ is a random tree with 
the same distribution as $T(p)$ conditioned on being infinite, 
we observe that for any 
$k=(k_1,\ldots,k_J)\in \{1,2,\ldots\}^J$
\begin{eqnarray} \label{eqn:Tinfty}
\lefteqn{\P(|T_\infty^{i_1}(p)|=\cdots=|T_\infty^{i_I}(p)|=\infty,
|T_\infty^{j_1}(p)|=k_1,\ldots,|T_\infty^{j_J}(p)|=k_J,
|T_\infty^{I+J+1}(p)|=0)}\nonumber\\
& & =\frac{(p\eta_\infty(p))^I\prod_{l=1}^Jp\eta_{k_l}(p)}{\eta_\infty(p)}
(1-p)=p(1-p)(p\eta_\infty(p))^{I-1}\prod_{l=1}^Jp\eta_{k_l}(p).
\end{eqnarray}
We now need to show that when using the construction of $\tilde{T}(p)$ we
get the analogous expression.
For any $i\in \cI$ let $j(i)$ be the smallest $j\in \cJ$ such that $j>i$
if such a $j$ exists. 
We get that
\begin{eqnarray} \label{eqn:long1}
\lefteqn{\P(|\tilde{T}^{i_1}(p)|=\cdots=|\tilde{T}^{i_I}(p)|=\infty,
|\tilde{T}^{j_1}(p)|=k_1,\ldots,|\tilde{T}^{j_J}(p)|=k_J,
|\tilde{T}^{I+J+1}(p)|=0)} \nonumber\\
& & =\sum_{l=1}^I 2^{-i_l} 
\P(L_{i_1}(p)=\cdots=L_{i_I}(p)=\infty,L_{j_1}(p)=k_1,\ldots,L_{j_J}(p)=k_J,
L_{I+J+1}(p)=0 | X=i_l) \nonumber\\
& & =\sum_{l=1}^I 2^{-i_l}(p\eta_\infty(p))^{l-1}
\left(\prod_{j_m<i_l}2p\eta_{k_m}(p)\right) \\
& & \ \ \ \ \times \P(L_{i_{l+1}}(p)=\cdots=L_{i_I}(p)=\infty,
           L_{j(i_l)}(p)=k_{j(i_l)},\ldots,L_{j_J}(p)=k_J,
L_{I+J+1}(p)=0 | X=i_l). \nonumber
\end{eqnarray}
The first equality simply divides into cases depending on the value of $X.$ 
The second equality uses (\ref{eqn:tildeT1}) and the first part 
of Lemma \ref{lemma:welldefined}. We proceed by considering (the possibly
empty) set of $j$ such that $i_l<j<i_{l+1}.$ 
Using that $2^{-i_l}\prod_{j_m<i_l}2=2^{-l}$, 
we then get that \eqref{eqn:long1} equals
\begin{eqnarray*}
\lefteqn{\sum_{l=1}^I 
2^{-l}(p\eta_\infty(p))^{l-1}\left(\prod_{j_m<i_l}p\eta_{k_m}(p)\right)
\left(\prod_{i_l<j_m<i_{l+1}}p\eta_{k_m}(p)\right)} \\
& & \ \ \ \ \times \P(L_{i_{l+1}}(p)=\cdots=L_{i_I}(p)=\infty,
            L_{j(i_{l+1})}(p)=k_{j(i_{l+1})},\ldots,L_{j_J}(p)=k_J,
           L_{I+J+1}(p)=0 | X=i_l)\\
 & & =\sum_{l=1}^I 2^{-l}(p\eta_\infty(p))^{l-1}
 \left(\prod_{j_m<i_{l+1}}p\eta_{k_m}(p)\right) p\eta_\infty(p)
 \frac{(2^l-1)p}{(2^l-2)p+1}\\
& & \ \ \ \ \times \P(L_{i_{l+2}}(p)=\cdots=L_{i_I}(p)=\infty,
           L_{j(i_{l+1})}(p)=k_{j(i_{l+1})},\ldots,L_{j_J}(p)=k_J,
           L_{I+J+1}(p)=0 | X=i_l),\\          
\end{eqnarray*}
where we in the second step  use Lemma \ref{lemma:welldefined}.
Iterating this procedure, and using \eqref{eqn:tildeT2}, we get that
\begin{eqnarray*} 
\lefteqn{\P(|\tilde{T}^{i_1}(p)|=\cdots=|\tilde{T}^{i_I}(p)|=\infty,|\tilde{T}^{j_1}(p)|=k_1,\ldots,|\tilde{T}^{j_J}(p)|=k_J,
|\tilde{T}^{I+J+1}(p)|=0)} \nonumber\\
& & =\cdots =\left(\prod_{m=1}^Jp\eta_{k_m}(p)\right)(p\eta_\infty(p))^{I-1}
\sum_{l=1}^I 2^{-l}\left(\prod_{m=l}^{I-1}
\frac{(2^m-1)p}{(2^m-2)p+1}\right)\frac{2^{I}p(1-p)}{(2^I-2)p+1},
\end{eqnarray*}
where 
\begin{equation} \label{eqn:equal1}
\prod_{m=l}^{I-1}\frac{(2^m-1)p}{(2^m-2)p+1}:=1,
\end{equation}
if $l=I.$   

Comparing this to (\ref{eqn:Tinfty}), we need to prove that for every $I,$
\begin{equation} \label{eqn:induction}
\sum_{l=1}^I 2^{-l}\left(\prod_{m=l}^{I-1}\frac{(2^m-1)p}{(2^m-2)p+1}\right)\frac{2^{I}}{(2^I-2)p+1}=1,
\end{equation}
and this we do by induction. First, we note that this trivially holds for $I=1.$ Assume therefore that it holds for
$I\geq 1,$ and observe that by \eqref{eqn:equal1},
\begin{eqnarray*}
\lefteqn{\sum_{l=1}^{I+1} 2^{-l}\left(\prod_{m=l}^{I}\frac{(2^m-1)p}{(2^m-2)p+1}\right)\frac{2^{I+1}}{(2^{I+1}-2)p+1}}\\
& & =\frac{2^{I+1}}{(2^{I+1}-2)p+1}\left(2^{-(I+1)}+\sum_{l=1}^{I} 2^{-l}\left(\prod_{m=l}^{I}\frac{(2^m-1)p}{(2^m-2)p+1}\right)\right)\\
& & =\frac{2^{I+1}}{(2^{I+1}-2)p+1}\left(2^{-(I+1)}+\sum_{l=1}^{I} 2^{-l}\left(\prod_{m=l}^{I-1}\frac{(2^m-1)p}{(2^m-2)p+1}\right)
\frac{2^{I}}{(2^I-2)p+1}\frac{(2^I-1)p}{2^I}\right)\\
& & =\frac{2^{I+1}}{(2^{I+1}-2)p+1}\left(2^{-(I+1)}+\frac{(2^I-1)p}{2^I}\right)
=\frac{1+2(2^I-1)p}{(2^{I+1}-2)p+1}=1,
\end{eqnarray*}
where we use the induction assumption in the third equality.
\fbox{}\\

In order to prove Theorem \ref{thm:main} using our construction of 
$\tilde{T}(p),$ we will need the monotonicity properties stated in our
next lemma, which consists of three parts.
\begin{lemma} \label{lemma:mon}
Consider the functions 
\begin{equation} \label{eqn:fcn1}
p\eta_k(p),
\end{equation}
and
\begin{equation} \label{eqn:fcn2}
\frac{2^np(1-p)}{(2^n-2)p+1}.
\end{equation}
\begin{enumerate}
\item For any $1\leq k<\infty,$ the function of \eqref{eqn:fcn1} 
is non-increasing in $p$ for $p\geq 1/2.$

\item For any $n\geq 1,$ the function of \eqref{eqn:fcn2} 
is non-increasing in $p$ for $p\geq 1/2.$

\item For any $1/2\leq p < 1,$ 
the function of \eqref{eqn:fcn2} is non-increasing in $n$ for 
$n \geq 1.$ 

\end{enumerate}
\end{lemma}
{\bf Proof.} We have that $p\eta_k(p)=c_{k}p^{k}(1-p)^k,$ which is clearly
non-increasing in $p$ for $p\geq 1/2.$ Furthermore, 
since $p(1-p)$ and $1/((2^n-2)p+1)$ are non-increasing in $p$ for $p\geq 1/2$ 
and any $n\geq 1,$ the first two parts of the statement follows.

The third part follows by observing that the function of \eqref{eqn:fcn2} 
is non-increasing in $n$ iff
\[
\left((2^{n+1}-2)p+1\right)2^{n}\geq 2^{n+1}\left((2^{n}-2)p+1\right),
\]
which simplifies to $1-2p \geq 2-4p$ and holds for all $p \in [1/2,1].$
\fbox{}\\

\section{Proof of Theorem \ref{thm:main}} \label{sec:mainproof}

In order to prove Theorem \ref{thm:main}, we will start by 
proving it in a special case, and then use this to prove the full statement.
\begin{theorem} \label{thm:speccase}
For any $p>1/2,$ there exists a coupling of $T_\infty(1/2)$ and 
$T_\infty(p)$ (where $T_\infty(1/2)\sim \cT_\infty(1/2)$ and
$T_\infty(p)\sim \cT_\infty(p)$) such that 
\[
\P(T_\infty(1/2) \subset T_\infty(p))=1.
\]
\end{theorem}
{\bf Remark:} The proof of Theorem \ref{thm:main} will be very similar, 
and therefore we will only address the necessary changes.

\medskip

\noindent
{\bf Proof of Theorem \ref{thm:speccase}.}
In order to facilitate the display of formulas, we will use the notation 
$p_1=1/2$ and $p_2=p.$

We will construct a sequence of pairs of trees 
$(\tilde{T}_l(p_1),\tilde{T}_l(p_2))_{l\geq 1}$ such that 
$\tilde{T}_l(p_i) \sim \cT_\infty(p_i)$ for every $l\geq 1,$ and
\begin{equation} \label{eqn:lorder}
\P(\tilde{T}_l(p_1)\cap(\cup_{n=0}^l \V_n)\subset 
\tilde{T}_l(p_2)\cap(\cup_{n=0}^l \V_n))=1,
\end{equation} 
for every $l\geq 1.$
That is, $(\tilde{T}_l(p_1),\tilde{T}_l(p_2))$ will be ordered 
up to distance $l$ from the root. From this, the statement will easily follow.

We start by proving \eqref{eqn:lorder} for $l=1.$
We will use the procedure of Section \ref{sec:prel} to generate our
trees, but we will do it simultaneously for $p_1$ and $p_2.$
Consider therefore the following random variables that we will use
in our construction of $(\tilde{T}_1(p_1),\tilde{T}_1(p_2)).$
\begin{enumerate}
\item $X$ is such that $\P(X=l)=2^{-l}$ for $l\geq 1.$ It is also 
independent of all other random variables listed.

\item $(U_i)_{i\geq 1}$ is an i.i.d. collection of $U[0,1]$ random variables
which is also independent of all other random variables
listed.

\item For every $m\geq 1,$ $(T_{k,m})_{k\geq 1}$ is a collection of 
random variables such that $T_{k,m} \sim \cT_k$ and 
$\P(T_{1,m}\subset T_{2,m}\subset \cdots)=1,$ which is possible
by Theorem \ref{thm:LW}.
Furthermore, we take  $(T_{k,m})_{k\geq 1}$ to be independent for
different $m,$ and independent from all other random variables
listed. We also let $T_{\infty,m}(p_1)=\bigcup_{k=1}^\infty T_{k,m}.$

\item Finally, $(T_{\infty,m}(p_2))_{m\geq 1}$
is an i.i.d. collection such that $T_{\infty,m}(p_2) \sim \cT_\infty(p_2)$ 
for every $m,$ and $(T_{\infty,m}(p_2))_{m\geq 1}$ is independent
from all other random variables listed.

\end{enumerate}

For any $m<X$ and $1\leq k<\infty,$ we let $L_m(p_i)=k$ iff
\[
\sum_{l=1}^{k-1}2p_i \eta_l(p_i) \leq U_m \leq \sum_{l=1}^{k}2p_i \eta_l(p_i)
\]
and otherwise we let $L_m(p_i)=\infty.$ Thus, we 
see that $L_m(p_i)$ is chosen with probabilities as in \eqref{eqn:tildeT1}
and in Lemma \ref{lemma:welldefined}.
Furthermore, by the first part of Lemma \ref{lemma:mon} we conclude
that $L_m(p_1)\leq L_m(p_2)$ for every such $m.$ We also let 
$L_X(p_1)= L_X(p_2)=\infty.$ Recall that 
$n_\infty(p_i,m)=|\{j\leq m-1:L_j(p_i)=\infty\}|$
and note that by our construction,
\begin{equation} \label{eqn:infineq}
n_\infty(p_1,X)\leq n_\infty(p_2,X).
\end{equation}
We proceed as follows.
For $m=X+1,$ we let $L_m(p_i)=0$ iff
\[
U_m < \frac{2^{n_\infty(p_i,m)}p_i(1-p_i)}{(2^{n_\infty(p_i,m)}-2)p_i+1},
\]
and for $1\leq k < \infty,$ we let $L_m(p_i)=k$ iff 
\[
\frac{2^{n_\infty(p_i,m)}p_i(1-p_i)}{(2^{n_\infty(p_i,m)}-2)p_i+1}+
\sum_{l=1}^{k-1}p_i\eta_l(p_i)\leq U_m <
\frac{2^{n_\infty(p_i,m)}p_i(1-p_i)}{(2^{n_\infty(p_i,m)}-2)p_i+1}+
\sum_{l=1}^{k}p_i\eta_l(p_i),
\]
and otherwise we let $L_m(p_i)=\infty.$
By Lemma \ref{lemma:mon}, $L_m(p_1)\leq L_m(p_2)$ and in particular
it is possible for $L_m(p_1)=0$ while $L_m(p_2)>0$.
This uses all three parts of that lemma.
As a consequence, we get that 
$n_\infty(p_1,X+1)\leq n_\infty(p_2,X+1).$

Continue for $m=X+2$ etc in the natural way,
and define $m_0(p_i)=\min\{m:L_m(p_i)=0\}.$
As above, Lemma \ref{lemma:mon}, shows that
$L_m(p_1)\leq L_m(p_2)$ for every $m$. 
Define 
\begin{equation} \label{eqn:tildeconstruct}
V(\tilde{T}_1(p_i)):=\{o\}\bigcup_{m=1}^{m_0(p_i)-1}
\bigcup_{v\in T_{L_m(p_i),m}}\{(m,v)\},
\end{equation}
where again, we abuse notation by writing $T_{L_m(p_i),m}$
instead of $T_{L_m(p_i),m}(p_i)$ when $L_m(p_i)=\infty.$

Thus, $\tilde{T}_1(p_i)$ is defined through the procedure of 
Section \ref{sec:prel},
and so by Lemmas \ref{lemma:tildeTdistr1} and \ref{lemma:tildeTdistr2} 
we see that 
$\tilde{T}(p_i)\sim \cT_\infty(p_i)$ for $i=1,2.$
By construction, we have that 
\[
\P(\tilde{T}_1(p_1)\cap(\{o\}\cup \V_1)\subset 
\tilde{T}_1(p_2)\cap(\{o\}\cup \V_1))=1.
\]

We now proceed inductively, and so assume that \eqref{eqn:lorder} holds for 
some $l\geq 1.$ We repeat the above procedure, but using 3' below instead of
points 3 and 4.

\begin{enumerate}
\item[3']   For every $m\geq 1,$ $(T_{k,m})_{k\geq 1}$ is a collection of 
random variables such that $T_{k,m} \sim \cT_k$ and 
$\P(T_{1,m}\subset T_{2,m}\subset \cdots)=1,$ which is possible
by Theorem \ref{thm:LW}.
We take  $(T_{k,m})_{k\geq 1}$ to be independent for
different $m,$ and also independent from all other random variables
listed, and let $T_{\infty,m}(p_1)=\bigcup_{k=1}^\infty T_{k,m}.$
Furthermore, we let $(T_{\infty,m}(p_2))_{m \geq 1}$ be such that 
$T_{\infty,m}(p_2) \sim \cT_{\infty}(p_2)$ and coupled with 
$T_{\infty,m}(p_1)$ so that 
\[
\P(T_{\infty,m}(p_1)\cap(\cup_{n=0}^l \V_n)\subset 
T_{\infty,m}(p_2)\cap(\cup_{n=0}^l \V_n))=1.
\]
This is possible by the induction hypothesis.
\end{enumerate}

We then construct $(\tilde{T}_{l+1}(p_1),\tilde{T}_{l+1}(p_2))$ using 
1,2 and 3'.
To see why $(\tilde{T}_{l+1}(p_1),\tilde{T}_{l+1}(p_2))$
satisfies \eqref{eqn:lorder} for $l+1,$ it suffices to observe the following.
By the induction hypothesis, any pair of trees used in \eqref{eqn:tildeconstruct} 
are already ordered up to distance $l$ from their roots, and so by attaching
them to vertices at distance one from the root of $\tilde{T}_{l+1}(p_i),$ 
the new trees are ordered up to distance $l+1$ from their roots.

In order to conclude the argument, let $\gamma_l$ be a measure on 
$\{0,1\}^{\V}\times \{0,1\}^{\V}$ with marginal distributions 
$\cT_{\infty}(p_1)$ and $\cT_{\infty}(p_2)$ such that 
$\gamma_l(\xi(\cup_{n=0}^l \V_n)\leq \eta(\cup_{n=0}^l \V_n))=1.$
Here, we identify a tree $T$ and an element $\xi_T \in \{0,1\}^{\V}$
by letting $\xi_T(v)=1$ iff $v\in T,$ and the measure $\gamma_l$ exists
by the above construction.
Since $\{0,1\}^{\V}\times \{0,1\}^{\V}$
is compact, there exists a subsequential 
limiting measure $\gamma$ on $\{0,1\}^{\V}\times \{0,1\}^{\V}$ with marginal 
distributions $\cT_{\infty}(p_1)$ and $\cT_{\infty}(p_2)$ such that 
$\gamma(\xi(\V)\leq \eta(\V))
=\lim_l \gamma(\xi(\cup_{n=0}^l \V_n)\leq \eta(\cup_{n=0}^l \V_n))=1.$
By Strassen's theorem it follows that there exists random trees 
$T_\infty(p_1)\sim \cT_\infty(p_1)$ and $T_\infty(p_2)\sim \cT_\infty(p_2)$
such that $\P(T_\infty(p_1)\subset T_\infty(p_2))=1.$
\fbox{}\\

We are now ready to prove Theorem \ref{thm:main}. Since the proof is very similar 
to the proof of Theorem \ref{thm:speccase}, we will only address 
the minor changes.

\noindent
{\bf Proof of Theorem \ref{thm:main}.} 
We can assume that $p_1>1/2$ since the case $p_1=1/2$ has already been proved.
We change point 3 to

\begin{enumerate}
\item[3] For every $m\geq 1,$ 
$(T_{k,m})_{k\geq 1}$ and $T_{\infty,m}(p_1)$ is a collection of 
random variables such that $T_{k,m} \sim \cT_k,$ 
$T_{\infty,m}(p_1)\sim \cT_{\infty}(p_1)$ and coupled so that 
$\P(T_{1,m}\subset T_{2,m}\subset \cdots \subset T_{\infty}(p_1) )=1,$ 
which is possible by Theorems \ref{thm:LW} and \ref{thm:speccase}.
Furthermore, we take  $(T_{k,m})_{k\geq 1}$ and $T_{\infty,m}(p_1)$
to be independent for
different $m,$ and also independent from all other random variables
listed.

\end{enumerate}
Point 3' is changed accordingly.

\fbox{}\\

{\bf Acknowledgement} I would like to thank Russ Lyons for suggesting the 
problem and an anonymous referee 
for useful comments and suggestions.

\end{document}